\documentclass{elsart}
\usepackage{amssymb}
\newtheorem{lemma}{Lemma}[section]
\newtheorem{theorem}{Theorem}[section]

\begin{document}
\begin{frontmatter}
\title{On a mixed cubic-superlinear non radially symmetric Schr\"odinger system - Part II: Numerical solutions}
\author{Abdurahman F. Aljohani}
\ead{a.f.aljohani@ut.edu.sa}
\address{Department of Mathematics, Faculty of Science, University of Tabuk, Saudi Arabia.}
\author{Anouar Ben Mabrouk\corauthref{cor1}}
\ead{anouar.benmabrouk@fsm.rnu.tn; amabrouk@ut.edu.sa}
\corauth[cor1]{Corresponding author}
\address{Research Unit of Algebra, Number Theory and Nonlinear Analysis UR11ES50, Department of Mathematics, Faculty of Sciences, 5019 Monastir. Tunisia.\\
\&
Department of Mathematics, Higher Institute of Applied Mathematics and Informatics, Street of Assad Ibn Al-Fourat, Kairouan 3100, Tunisia.\\
\& Department of Mathematics, Faculty of Science, University of Tabuk, Saudi Arabia.}
\begin{abstract}
In this paper a nonlinear coupled Schrodinger system in the presence of mixed cubic and superlinear power laws is considered. A non standard numerical method is developed to approximate the solutions in higher dimensional case. The idea consists in transforming the continuous system into an algebraic quasi linear dynamical discrete one leading to generalized semi-linear operators. Next, the discrete algebraic system is studied for solvability, stability, convergence and stability. At the final step, numerical examples are provided to illustrate the efficiency of the theoretical results. 
\end{abstract}
\end{frontmatter}
\section{Introduction}
The present work is devoted to the numerical study of a coupled system of nonlinear Schr\"odinger equations characterized by a mixed nonlinearities. Focuses are made on the development of a non standard numerical method will be developed to study the numericall solutions of the original evolutive system by means of sophisticated algebraic operators such as the famous Lyapunov-Sylvester ones in a higher dimensional case. 

Denote for $\lambda$ and $p$ real numbers such that $\lambda>0$ and $p>1$,
$$
g(u,v)=|u|^{p-1}+\lambda|v|^{2}\;\mbox{and}\;f(u,v)=g(u,v)u.
$$
We consider in the first part the evolutive system
\begin{equation}\label{ContinuousProblem1}
\left\{\begin{array}{lll}
iu_{t}+\sigma_1\Delta u+g(u,v)u=0,\\
iv_{t}+\sigma_2\Delta v+g(v,u)v=0
\end{array}
\right.
\end{equation}
with the extra initial conditions
\begin{equation}\label{ContinuousInitialCondition}
\left\{\begin{array}{lll}
W(x,y,t_{0})=W_{0}(x,y)=(u_0(x,y),v_0(x,y))(x,y)\in\Omega\\ \mbox{and}\\
\frac{\partial W}{\partial t}(x,y,t_{0})=W_1(x,y)=(u_1(x,y),v_1(x,y)),\;(x,y)\in\Omega
\end{array}
\right.
\end{equation}
and boundary conditions
\begin{equation}\label{ContinuousBoundaryCondition}
\frac{\partial W}{\partial\eta}(x,y,t)=0,\quad((x,y),t)\in\partial\Omega\times(t_{0},+\infty)
\end{equation}
where $\Omega =[L_{0},L_1]\times[L_{0},L_1]$ is a rectangular domain in $\mathbb{R}^2$. \\
$u_{t}$ is the first order partial derivative in time, $u_{tt}$ is the second order partial derivative in time, $\Delta =\frac{\partial ^2}{\partial x^2}+\frac{\partial ^2}{\partial y^2}$ is the Laplace operator on $\mathbb{R}^2$. $\frac{\partial}{\partial\eta}$ is the outward normal derivative operator along the boundary $\partial \Omega $. $u_{0}$ and $u_1$ are real valued functions with $u_{0}$ and $u_1$ are $\mathcal{C}^2$ on $\overline{\Omega}$. $u$ and $v$ are the unknown candidates supposed to be $\mathcal{C}^{4}$ on $\overline{\Omega}$, $W=(u,v)$. $\sigma_i$, $i=1,2$ are real parameters such that $\sigma_i>0$.

We intend to apply generalized Lyapunov-Sylvester algebraic operators instead of transforming the two-dimensional discrete problem into block-tridiagonal form, to develop numerical solutions for the system (\ref{ContinuousProblem1})-(\ref{ContinuousBoundaryCondition}) by replacing time and space partial derivatives by finite-difference approximations. The used method is better as it leads to fast convergent and more accurate discrete algebraic systems. It permits also to somehow conserve the geometric presentation of the problem as we  solve in the same two-dimensional space and did not project the problem on one-dimensional grids. Relatively to computer architecture, the process of projecting on different spaces and next lifting to the original one may induce degradation of error estimates and slow algorithms.

The next section is concerned with the study of the numerical solutions of the system (\ref{ContinuousProblem1})-(\ref{ContinuousBoundaryCondition}). More precisely an introductory part is devoted to the introduction of the discretization method. Next, solvability of such a discrete system is proved in section 3. Section 4 is concerned with the consistency, stability and the convergence of the discrete Lyapunov-Sylvester problem obtained in section 3 by applying the truncation error for consistency, Lyapunov cretirion for stability and the Lax equivalence theorem for the convergence. Section 5 is devoted to the development of numerical examples. Performance of the discrete scheme is proved by means of error estimates as well as fast algorithms. The conclusion is finally subject of section 6.

\section{Discrete two-dimensional nonlinear NLS system}
The object of this section is develop a discretization scheme to approximate numerically the solution(s) of the evolutive (time-dependent) problem (\ref{ContinuousProblem1})-(\ref{ContinuousBoundaryCondition}). The proposed schme permits to transform problem (\ref{ContinuousProblem1})-(\ref{ContinuousBoundaryCondition}) into a discrete quasi-linear one which by the next will be studied for convergence, solvability and consistency. Consider a time step $l=\Delta t$ and a space one $h=\displaystyle\frac{L_1-L_0}{J+1}$. Next, denote for $n\in\mathbb{N}$ and $j,m\in\{0,...,J\}$
$$
t^n=t_0+nl\;,x_j=L_0+jh\quad\hbox{and}\quad y_m=L_0+mh
$$
so that the cube $[L_0,L_1]\times[L_0,L_1]$ is subdivided into cubes $C_{j,m}=[x_j,x_{j+1}]\times[y_m,y_{m+1}]$. For a function $z$ defined on the cube $[L_0,L_1]\times[L_0,L_1]$, we denote by small $z_{j,m}^n$ the net function $z(x_j,y_m,t^n)$ and capital $Z_{j,m}^n$ the numerical approximation. Consider next the discrete finite difference operators
$$
u_t=\displaystyle\frac{u^{n+1}-u^{n-1}}{2\ell},
$$
$$
u_x=\displaystyle\frac{\overline{u}^{n}_{j+1,m}-\overline{u}^{n}_{j-1,m}}{2h},\;\;\;
u_y=\displaystyle\frac{\overline{u}^{n}_{j,m+1}-\overline{u}^{n}_{j,m-1}}{2h},
$$
$$
\Delta u
=\displaystyle\frac{\overline{u}^{n}_{j+1,m}-2\overline{u}^{n}_{j,m}+\overline{u}^{n}_{j-1,m}}{h^2}
+\displaystyle\frac{\overline{u}^{n}_{j,m+1}-2\overline{u}^{n}_{j,m}+\overline{u}^{n}_{j,m-1}}{h^2},
$$
$$
\overline{u}^n=\mu_1u^{n+1}+\mu_2u^n+\mu_3u^{n-1},
$$
with $\mu_i\in(0,1)$, $i=1,2,3$ such that $\mu_1+\mu_2+\mu_3=1$. It is a barycentric calibration method that is applied firstly by Bratsos and his collaborators and which has been proved to be useful and efficient. In fact, it is always questionable to confirm what is the closest to the exact (unknown) value $u(x,y,t)$ on the grid $U^n_{j,m}$. This explains the use of the calibration proposed. (See \cite{Bezia-BenMabrouk-Betina1}, \cite{Bezia-BenMabrouk1}, \cite{Bratsos1}, \cite{Bratsos2}, \cite{Chteouietal1}).

Using the discrete operators introduced above, problem (\ref{ContinuousProblem1}) may be written on the discrete form
\begin{equation}\label{DiscreteProblem1}
\left\{
\begin{array}{lll}
\begin{array}{lll}
&&i\displaystyle\frac{u^{n+1}_{j,m}-u_{j,m}^{n-1}}{\ell}
+\sigma_1\displaystyle\frac{\overline{u}^{n+1}_{j+1,m}-2\overline{u}^{n+1}_{j,m}+\overline{u}_{j-1,m}^{n+1}}{h^2}\\
&&+\sigma_1\displaystyle\frac{\overline{u}^{n+1}_{j,m+1}-2\overline{u}^{n+1}_{j,m}+\overline{u}_{j,m-1}^{n+1}}{h^2}\\
&&+g(u^{n}_{j,m},v^{n}_{j,m})\overline{u}_{j,m}^{n}=0.
\end{array}\\
\\
\begin{array}{lll}
&&i\displaystyle\frac{v^{n+1}_{j,m}-v_{j,m}^{n-1}}{\ell}
+\sigma_2\displaystyle\frac{\overline{v}^{n+1}_{j+1,m}-2\overline{v}^{n+1}_{j,m}+\overline{v}_{j-1,m}^{n+1}}{h^2}\\
&&+\sigma_2\displaystyle\frac{\overline{v}^{n+1}_{j,m+1}-2\overline{v}^{n+1}_{j,m}+\overline{v}_{j,m-1}^{n+1}}{h^2}\\
&&+g(v^{n}_{j,m},u^{n}_{j,m})\overline{v}_{j,m}^{n}=0.
\end{array}
\end{array}\right.
\end{equation}
Denote next $\sigma=\displaystyle\frac{l}{h^2}$. The first equation in system (\ref{DiscreteProblem1}) may be written as
$$
\begin{array}{lll}
&&i(u_{j,m}^{n+1}-u_{j,m}^{n-1})+\sigma_1\sigma\left[\mu_1u_{j+1,m}^{n+1}+\mu_2u_{j+1,m}^{n}+\mu_3u_{j+1,m}^{n-1}\right.\\
&&-2\mu_1u_{j,m}^{n+1}-2\mu_2u_{j,m}^{n}-2\mu_3u_{j,m}^{n-1}
+\mu_1u_{j-1,m}^{n+1}+\mu_2u_{j-1,m}^{n}+\mu_3u_{j-1,m}^{n-1}\\
&&+\mu_1u_{j,m+1}^{n+1}+\mu_2u_{j,m+1}^{n}+\mu_3u_{j,m+1}^{n-1}\\
&&\left.-2\mu_1u_{j,m}^{n+1}-2\mu_2u_{j,m}^{n}-2\mu_3u_{j,m}^{n-1}+\mu_1u_{j,m-1}^{n+1}+\mu_2u_{j,m-1}^{n}+\mu_3u_{j,m-1}^{n-1}\right]\\
&&+g(u_{j,m}^n,v_{j,m}^n)\left(\mu_1u_{j,m}^{n+1}+\mu_2u_{j,m}^{n}+\mu_3u_{j,m}^{n-1}\right)=0.
\end{array}
$$
Denote now for $i,j=1,2,3$,
$$
\Gamma_{ij}=\sigma_i\mu_j\sigma
$$
and for $j,m=0,1,\dots,J$,
$$
\Gamma_{j,m}^{n,1}=\displaystyle\frac{1}{2}(g(j,m)^n-4\sigma_1\sigma),\;\;
$$
$$ 
\Lambda_{j,m}^{n}=\displaystyle\frac{1}{2}(i+2\mu_1\Gamma_{j,m}^{n,1})
$$
and
$$
\widetilde{\Lambda_{j,m}^{n}}=\displaystyle\frac{1}{2}(-i+2\mu_3\Gamma_{j,m}^{n,1}).
$$
We get
$$
\begin{array}{lll}
&&\Gamma_{11}u_{j-1,m}^{n+1}+\Lambda_{j,m}^{n}u_{j,m}^{n+1}+\Gamma_{11}u_{j+1,m}^{n+1}\\
&+&\Gamma_{11}u_{j,m-1}^{n+1}+\Lambda_{j,m}^{n}u_{j,m}^{n+1}+\Gamma_{11}u_{j,m+1}^{n+1}\\
&+&\Gamma_{12}u_{j-1,m}^{n}+\mu_2\Gamma_{j,m}^{n,1}u_{j,m}^{n}+\Gamma_{12}u_{j+1,m}^{n}\\
&+&\Gamma_{12}u_{j,m-1}^{n}+\mu_2\Gamma_{j,m}^{n,1}u_{j,m}^{n}+\Gamma_{12}u_{j,m+1}^{n}\\
&+&\Gamma_{13}u_{j-1,m}^{n-1}+\widetilde{\Lambda}_{j,m}^{n}u_{j,m}^{n-1}+\Gamma_{13}u_{j+1,m}^{n-1}\\
&+&\Gamma_{13}u_{j,m-1}^{n-1}+\widetilde{\Lambda}_{j,m}^{n}u_{j,m}^{n-1}+\Gamma_{13}u_{j,m+1}^{n-1}=0.
\end{array}
$$
Exploiting the boundary conditions, the last equation may be written in a matrix-vector form
\begin{equation}\label{LyapEq1}
A_1^nU^{n+1}+U^{n+1}A_1^n+A_2^nU^{n}+U^{n}A_2^n+A_3^nU^{n-1}+U^{n-1}A_3^n=0,
\end{equation}
where $U^n=(u_{j,m}^n)$ and $V^n=(v_{j,m}^n)$ are the unknown solutions and for $i=1,2,3$ $A_i^n$ are the matrices given by
$$
A_1^n(0,1)=A_1^n(J,J-1)=2\Gamma_{11},
$$
$$
A_1^n(j,j)=\Lambda_{j,m}^{n}\;,\;\;0\leq j\leq J,
$$
$$
A_1^n(j,j+1)=A_1^n(j,j-1)=\Gamma_{11}\;,\;\;,\;\;1\leq j\leq J-1,
$$
$$
A_2^n(0,1)=A_2^n(J,J-1)=2\Gamma_{12},
$$
$$
A_2^n(j,j)=\mu_2\Gamma_{j,m}^{n,1}\;,\;\;0\leq j\leq J,
$$
$$
A_2^n(j,j+1)=A_2^n(j,j-1)=\Gamma_{12}\;,\;\;,\;\;1\leq j\leq J-1,
$$
$$
A_3^n(0,1)=A_3^n(J,J-1)=2\Gamma_{13},
$$
$$
A_3^n(j,j)=\widetilde{\Lambda_{j,m}^{n}}\;,\;\;0\leq j\leq J,
$$
$$
A_3^n(j,j+1)=A_3^n(j,j-1)=\Gamma_{13}\;,\;\;,\;\;1\leq j\leq J-1,
$$
Similarly, the second eqaution in (\ref{DiscreteProblem1}) may be written on the form
$$
\begin{array}{lll}
&&i(v_{j,m}^{n+1}-v_{j,m}^{n-1})+a_2\sigma\left[\mu_1v_{j+1,m}^{n+1}+\mu_2v_{j+1,m}^{n}+\mu_3v_{j+1,m}^{n-1}\right.\\
&&-2\mu_1v_{j,m}^{n+1}-2\mu_2v_{j,m}^{n}-2\mu_3v_{j,m}^{n-1}
+\mu_1v_{j-1,m}^{n+1}+\mu_2v_{j-1,m}^{n}+\mu_3v_{j-1,m}^{n-1}\\
&&+\mu_1v_{j,m+1}^{n+1}+\mu_2v_{j,m+1}^{n}+\mu_3v_{j,m+1}^{n-1}\\
&&\left.-2\mu_1v_{j,m}^{n+1}-2\mu_2v_{j,m}^{n}-2\mu_3v_{j,m}^{n-1}+\mu_1v_{j,m-1}^{n+1}+\mu_2v_{j,m-1}^{n}+\mu_3v_{j,m-1}^{n-1}\right]\\
&&+g(v_{j,m}^n,u_{j,m}^n)\left(\mu_1v_{j,m}^{n+1}+\mu_2v_{j,m}^{n}+\mu_3v_{j,m}^{n-1}\right)=0.
\end{array}
$$
Here also denote similarly
$$
\Gamma_{j,m}^{n,2}=\displaystyle\frac{1}{2}(g(m,j)^n-4\sigma_2\sigma),
$$
$$
\Theta_{j,m}^{n}=\displaystyle\frac{1}{2}(i+2\mu_1\Gamma_{j,m}^{n,2})
$$
nd
$$
\widetilde{\Theta_{j,m}^{n}}=\displaystyle\frac{1}{2}(-i+2\mu_3\Gamma_{j,m}^{n,2}).
$$
We get
$$
\begin{array}{lll}
&&\Gamma_{21}v_{j-1,m}^{n+1}+\Theta_{j,m}^nv_{j,m}^{n+1}+\Gamma_{21}v_{j+1,m}^{n+1}\\
&+&\Gamma_{21}v_{j,m-1}^{n+1}+\Theta_{j,m}^nv_{j,m}^{n+1}+\Gamma_{21}v_{j,m+1}^{n+1}\\
&+&\Gamma_{22}v_{j-1,m}^{n}+\mu_2\Gamma_{j,m}^{n,2}v_{j,m}^{n}+\Gamma_{22}v_{j+1,m}^{n}\\
&+&\Gamma_{22}v_{j,m-1}^{n}+\mu_2\Gamma_{j,m}^{n,2}v_{j,m}^{n}+\Gamma_{22}v_{j,m+1}^{n}\\
&+&\Gamma_{23}v_{j-1,m}^{n-1}+\widetilde{\Theta}_{j,m}^{n}v_{j,m}^{n-1}+\Gamma_{23}v_{j+1,m}^{n-1}\\
&+&\Gamma_{23}v_{j,m-1}^{n-1}+\widetilde{\Theta}_{j,m}^{n}v_{j,m}^{n-1}+\Gamma_{23}v_{j,m+1}^{n-1}=0.
\end{array}
$$
Exploiting the boundary conditions, the last equation may be written in a matrix-vector form
\begin{equation}\label{LyapEq2}
B_1^nV^{n+1}+V^{n+1}B_1^n+B_2^nV^{n}+V^{n}B_2^n+B_3^nV^{n-1}+V^{n-1}B_3^n=0.
\end{equation}
where for $i=1,2,3$, $B_i^n$ are the matrices given by
$$
B_1^n(0,1)=B_1^n(J,J-1)=2\Gamma_{21},
$$
$$
B_1^n(j,j)=\Theta_{j,m}^{n}\;,\;\;0\leq j\leq J,
$$
$$
B_1^n(j,j+1)=B_1^n(j,j-1)=\Gamma_{21}\;,\;\;,\;\;1\leq j\leq J-1,
$$
$$
B_2^n(0,1)=B_2^n(J,J-1)=2\Gamma_{22},
$$
$$
B_2^n(j,j)=\mu_2\Gamma_{j,m}^{n,2}\;,\;\;0\leq j\leq J,
$$
$$
B_2^n(j,j+1)=B_2^n(j,j-1)=\Gamma_{22}\;,\;\;,\;\;1\leq j\leq J-1,
$$
$$
B_3^n(0,1)=B_3^n(J,J-1)=2\Gamma_{23},
$$
$$
B_3^n(j,j)=\widetilde{\Theta_{j,m}^{n}}\;,\;\;0\leq j\leq J,
$$
$$
B_3^n(j,j+1)=B_3^n(j,j-1)=\Gamma_{23}\;,\;\;,\;\;1\leq j\leq J-1,
$$
\section{Solvability of the discrete problem}
Usually, discrete schemes used for numerical solutions of PDEs are transformed to algebraic equations on the form $AU^{n+1}=F(U^,U^{n-1},\dots,U^0)$ where $A$ is matrix or generally a linear operator. Next, the problem becomes whether this operator is invertible or not. The most known methods are based on eigenvalues/eignevectors computation of such operators. See \cite{Benmabrouk1}, \cite{Benmabrouk-Benmohamed-Omrani}, \cite{Bratsos1}, \cite{Bratsos2}, In the presentt work, we will not apply such procedure, but we develop different arguments based on the invertibility of Lyapunov-Sylvester operators as in \cite{Bezia-BenMabrouk-Betina1}, \cite{Bezia-BenMabrouk1} and \cite{Chteouietal1}. The first main result in this part is stated as follows.
\begin{theorem}\label{theorem1}
The system (\ref{LyapEq1})-(\ref{LyapEq2}) is uniquely solvable whenever the solutions $W^0=(U^0,V^0)$ and $W^1=(U^1,V^1)$ are known.
\end{theorem}
The proof reposes on the inverse of Lyapunov-Syslvester operators. Consider the endomorphism $\Phi$ defined by
\begin{equation}\label{LyapSylOperator}
\Phi_{l,h}^n(X,Y)=(\,A_1^nX+XA_1^n\,,\,B_1^nY+YB_1^n\,),
\end{equation}
To prove Theorem \ref{theorem1}, we need some preliminary results.
\begin{lemma}\label{LyapunovStabilityLemma}
	$\mathcal{P}_n$: The solution $(U^n,V^n) $ is bounded independently of $n$ whenever the initial solution $(U^0,V^0)$ is bounded.
\end{lemma}
\textbf{Proof.}
	Writing the initial condition in the discrete form we get
	\begin{equation}
	W^{2}=W^{0}+2lW^1.
	\end{equation}
	For $n=1$, this yields that $W^2$ is bounded. So assume next that
	\begin{equation}\label{hyporecurrence}
	W^k\mbox{ is bounded independently of }k\,;\;k=0,1,\dots,n.
	\end{equation}
	We shall show that $W^{n+1}$ is bounded independently of $n$. We already know from (\ref{LyapEq1}) that
	\begin{equation}\label{LyapEq1bis}
	\Phi_{l,h}^n(U^{n+1})=-\bigl(A_2^nU^{n}+U^{n}A_2^n+A_3^nU^{n-1}+U^{n-1}A_3^n\bigr),
	\end{equation}
	where $\Phi_{l,h}^n$ is the operator defined on the space of $(J+1,J+1)$-matrices $\mathcal{M}_{J+1}(\mathbb{C})$ by
	$$
	\Phi_{l,h}^n(X,Y)=A_1^nX+XA_1^n.
	$$
	From the recurrence hypothesis (\ref{hyporecurrence}), the matrice $A_1^n$ is bounded uniformly independently of $n$. Consequently, whenever $l=o(h^2)$ and $l,h\longrightarrow0$, we get
	\begin{equation}\label{limitofPhilhn}
	\Phi_{l,h}^n{\rightarrow}\,iId\;\mbox{as}\;{l,h\longrightarrow0}
	\end{equation}
	uniformly on $n$. As a consequence, there exists a constant $C=C(l,h)>0$, for which
	\begin{equation}\label{limitofPhilhnbis}
	\|\Phi_{l,h}^n\|\geq C,
	\end{equation}
	for $(l,h)$ small enough. It follows from (\ref{LyapEq1bis}) that
	\begin{equation}\label{hyprecurrencen+1}
	C\|U^{n+1}\|\leq\|A_2^nU^{n}+U^{n}A_2^n+A_3^nU^{n-1}+U^{n-1}A_3^n\|.
	\end{equation}
	The right hand term is bounded independently of $n$ from the recurrence hypothesis. As a result, $U^{n+1}$ is also bounded independently of $n$.\\
	The same result may be proved for $V^{n+1}$ by using equation (\ref{LyapEq2}) and the operator
	$$
	\widetilde{\Phi}_{l,h}^n(X,Y)=B_1^nX+XB_1^n.
	$$

Next, we apply the following result.
\begin{lemma} \label{LemmeInversion}
	Let $E$ be a finite dimensional ($\mathbb{R}$ or $\mathbb{C}$) vector space and $(\Phi_n)_n$ be a sequence of endomorphisms converging uniformly to an invertible endomorphism $\Phi$. Then, there exists $n_{0}$ such that, for any $n\geq\,n_{0}$, the endomorphism $\Phi_n$ is invertible.
\end{lemma}
\textbf{Proof of Theorem \ref{theorem1}.}
It follows from the arguments of Lemma \ref{LyapunovStabilityLemma} and Lemma \ref{LemmeInversion} that the operator defined by
$$
\Phi_{l,h}^n(X,Y)=(\Phi_{l,h}^n(X),\widetilde{\Phi}_{l,h}^n(Y))
$$
is an endomorphism for $l,h$ small enough. Which gives the desired result.
\section{Consistency, stability and convergence of the discrete method}
Recall firstly that the consistency of the numerical scheme is always done by evaluating the local truncation error arising from the discrete and the continuous problem. In the present case, we have the following lemma.
\begin{lemma}
	\begin{itemize}
		\item Whenever $\mu_1=\mu_3$, the discrete scheme is consistent with order $o(l^2+h^2)$.
		\item Whenever $\mu_1\not=\mu_3$, the discrete scheme is consistent with order $o(l+h^2)$.
	\end{itemize}
\end{lemma}
\textbf{Proof.} Applying Taylor's expansion in the discrete equations raised in section 2, we get the following truncation principal part for the first equation in system (\ref{ContinuousProblem1})
	$$
	\begin{array}{lll}
	\mathcal{L}_{u,v}^1(x,y,t)
	&=&(\mu_1-\mu3)\sigma_1\displaystyle\frac{\partial}{\partial t}(\Delta u)l+\displaystyle\frac{\mu_1+\mu3}{2}\sigma_1\displaystyle\frac{\partial^2}{\partial t^2}(\Delta u)l^2\\
	&+&\displaystyle\frac{\sigma_1}{12}(\Delta^2u)h^2+o(l^2+h^2).
	\end{array}
	$$
	and for the second equation, we get
	$$
	\begin{array}{lll}
	\mathcal{L}_{u,v}^2(x,y,t)
	&=&(\mu_1-\mu3)\sigma_2\displaystyle\frac{\partial}{\partial t}(\Delta v)l+\displaystyle\frac{\mu_1+\mu3}{2}\sigma_2\displaystyle\frac{\partial^2}{\partial t^2}(\Delta v)l^2\\
	&+&\displaystyle\frac{\sigma_2}{12}(\Delta^2v)h^2+o(l^2+h^2).
	\end{array}
	$$
	where $\Delta_2=\displaystyle\frac{\partial^4}{\partial x^4}+\frac{\partial^4}{\partial y^4}$, $W=(u,v)$. Hence, the Lemma is proved.

Nex, the stability of the discrete scheme will be examined using the Lyapunov criterion of stability. Recall that a dynamical system $\mathcal{L}(u_{n+1},u_{n},u_{n-1},\dots)=0$ is stable in the sense of Lyapunov iff for any bounded initial value, the solution $u_n$ ramains bounded for all $n\geq0$. In the present case, we have the following result.
\begin{lemma}
	The discret system () is stabel in the sense of Lyapunov stability.
\end{lemma}
We already proved in Lemma \ref{LyapunovStabilityLemma} the property $\mathcal{P}_n$ affirming that the solution $W^n$ is bounded independently of $n$ whenever the initial solution $(W^0,W^1$ is bounded.

Now, it remains finally to check the convergence of the discrete scheme. This is done by a direct application of the following well-known result \cite{Lax-Richtmyer}.
\begin{theorem}
	(\textbf{Lax Equivalence Theorem}). For a consistent finite difference scheme, stability is equivalent to convergence.
\end{theorem}
\begin{lemma} \label{laxequivresult}
	As the numerical scheme is consistent and stable, it is then convergent.
\end{lemma}
\section{Numerical implementation}
We present in this section some illustrative examples in order to validate the methods and the results just described above. Recall that nonlinear Schr\"odinger equation plays an important role in the modeling of many phenomena. We mention as examples the models of Bose-Einstein condensation and the stabilized solitons. In the latter case, the nonlinear Schr\"odinger equation gives rise to soliton solutions in which the explicit expression can be well defined. For example, in the case of nonlinear cubic Schr\"odinger equation, $u$ is given by
$$
u(x,t)=K_u\exp\Bigl(i\bigl(\displaystyle\frac{1}{2}cx-\theta t+\varphi\bigr)\Bigr)sech\Bigl(\sqrt{a}(x-ct)+\phi\Bigr)
$$
where $a$, $q_s$, $c$, $\theta$, $\varphi$ and $\phi$ are some appropriate constants. For $t$ fixed, this function decays exponentially as $|x|\rightarrow\infty$. It is a soliton-type disturbance which travels with speed $c$ and with $a$-governed amplitude. For backgrounds on such a subject, the readers may refer to \cite{Bratsos1}, \cite{Bratsos2} and \cite{Lamb}.

Recall also that soliton type particles are always travelling along the whole real line, but as it is said with en exponential decay at the boundaries. So, to compute a solution, we need first to make some additional artificial hypothesis affirming that for some compact support $[L_0,L_1]$, we have $u(L_0,t)=u(L_1,t)=0$ for all $t$. Such hypothesis is not exact in general. However, many solutions have fast decay at infinity such as solitons. So one can reasonably use such it.

To measure the closeness of the numerical solution and the exact one, the error is evaluated via an $L_2$ matrix norm
$$
\|X\|_2=\Big(\sum_{i,j=0}^{J}|X_{ij}|^2\Big)^{\frac{1}{2}}
$$
for a matrix $X=(X_{ij})\in\mathcal{M}_{J+1}(\mathbb{C})$. Denote $u^n$ the net function $u(x,y,t^n)$ and $U^n$ the numerical solution. We propose to compute the discrete error
\begin{equation}
\mathrm{Er}=\max_n\|U^n-u^n\|_2  \label{Er}
\end{equation}
on the grid $(x_{i},y_{j})$, $0\leq\,i,j\leq J$ and the relative error between the exact solution and the numerical one as
\begin{equation}
\mathrm{Relative\,Er}=\max_n\frac{\|U^n-u^n\|_2}{\|u^n\|_2}  \label{Errelative}
\end{equation}
In the present paper, we consider the phenomena of propagation and interaction of solitons. We fix the problem parameters $p$, $\lambda$, $\sigma_1$ and $\sigma_2$ to be
$$
p=\displaystyle\frac{5}{2}\,,\;\;\lambda=\sigma_1=\sigma_2=1.
$$
\subsubsection{Simultaneous propagation of two solitons}
In this subsection, we illustrate numerical solutions of two solitons propagating simultaneously. For simplicity, denote
\begin{equation}\label{Soliton1-2-expression}
\left\{
\begin{array}{lll}
A(x,y,t)=\omega t-\displaystyle\frac{c}{2}x+\displaystyle\frac{c}{2}y+\varphi_v\,,\\
B(x,y,t)=sech(\sqrt{a}(x-y-ct)+\phi_u)\,,\\
C(x,y,t)=\omega t+\displaystyle\frac{c}{2}x-\displaystyle\frac{c}{2}y+\varphi_v\,,\\
D(x,y,t)=sech(\sqrt{a}(y-x-ct)+\phi_u)\,.
\end{array}
\right.
\end{equation}
The first is governed by the exact solution $u$ given by
\begin{equation}\label{Soliton1}
u(x,y,t)=K_u\exp\Bigl(iA(x,y,t)\Bigr)\Biggl[sech\Bigl(B(x,y,t)\Bigr)\Biggr]^{\frac{4}{3}}.
\end{equation}
It consists of a soliton traveling in the direction of the vector $(1,-1)$. The second soliton is governed by the exact solution
\begin{equation}\label{Soliton2}
v(x,y,t)=K_v\exp\Bigl(iC(x,y,t)\Bigr)\Biggl[sech\Bigl(D(x,y,t)\Bigr)\Biggr]^{\frac{4}{3}}.
\end{equation}
It consists of a soliton traveling in the direction of the same vector $(1,-1)$ but in the opposite direction of the first one, which guarantee in soliton physics the phenomenon of interaction. Denote also for any function $f=f(x,y,t)$
$$
T_f(x,y,t)=\tanh(f(x,y,t)).
$$
Technical calculus yield that
\begin{equation}\label{Soliton1-2-constants}
K_u=\Biggl(\displaystyle\frac{32a}{9}\Biggr)
^{\frac{2}{3}}\,,K_v=\Biggl(\displaystyle\frac{56a}{9}\Biggr)^{\frac{2}{3}}\,\mbox{and}\,\omega+\displaystyle\frac{c^2}{2}=\displaystyle\frac{36a}{3}.
\end{equation}
The computations are done for $-80\leq x,y\leq 100$ with different space and time steps as shown in the figures and tables corresponding. We fix also soliton parameters $a=0,01$, $c=0,1$ and the phase parameters $\varphi=\phi=0$. Therefore, the soliton pair $(u,v)$ is the exact solution of the inhomogeneous problem
\begin{equation}\label{InhpmogeneousContinuousProblem1}
\left\{
\begin{array}{lll}
iu_{t}+\Delta u+|u|^{p-1}u+|v|^2u=G_1(x,y,t),\\
iv_{t}+\Delta v+|v|^{p-1}v+|u|^2v=G_2(x,y,t)
\end{array}
\right.
\end{equation}
where $G_1$ and $G_2$ are explicited respectively by
$$
G_1(x,y,t)=\biggl[K_v^2\Bigl(1-T_D^2(x,y,t)\Bigr)^{\frac{4}{3}}-8aT_B^2(x,y,t)+i4c\sqrt{a}T_B(x,y,t)\biggr]u(x,y,t)
$$
and
$$
G_2(x,y,t)=\biggl[K_u^2\Bigl(1-T_B^2(x,y,t)\Bigr)^{\frac{4}{3}}+i\displaystyle\frac{4c\sqrt{a}}{3}T_D(x,y,t)\biggr]v(x,y,t).
$$
\subsubsection{Simultaneous propagation of x-soliton/y-soliton}
In this subsection, we try to illustrate the phenomenon of propagation of one first soliton propagating on the $x$-direction and an other one on the $y$-direction, which are already governed by the coupled system (\ref{ContinuousProblem1}) and the possible interaction between them. As previously, we denote for simplicity, 
\begin{equation}\label{Soliton1-2-expression}
\left\{\begin{array}{lll}
A(x,t)=\omega t+\displaystyle\frac{c}{2}x+\varphi_v\,,\\
B(x,t)=sech(\sqrt{a}(x-ct)+\phi_u)\,,\\
C(y,t)=\omega t-\displaystyle\frac{c}{2}y+\varphi_v\,,\\
D(y,t)=sech(\sqrt{a}(y-ct)+\phi_u)\,.
\end{array}
\right.
\end{equation}
The first is governed by the exact solution $u$ given by
\begin{equation}\label{Soliton1}
u(x,y,t)=u(x,t)=K\exp\Bigl(iA(x,t)\Bigr)\Biggl[sech\Bigl(B(x,t)\Bigr)\Biggr]^{\frac{4}{3}}.
\end{equation}
The second soliton is governed by the exact solution
\begin{equation}\label{Soliton2}
v(x,y,t)=v(y,t)=K\exp\Bigl(iC(y,t)\Bigr)\Biggl[sech\Bigl(D(y,t)\Bigr)\Biggr]^{\frac{4}{3}}.
\end{equation}
Denote also for any function $f=f(x,y,t)$
$$
T_f(x,y,t)=\tanh(f(x,y,t)).
$$
Technical calculus yield that
\begin{equation}\label{Soliton1-2-constants}
K=\Biggl(\displaystyle\frac{28a}{9}\Biggr)^{\frac{2}{3}}\,\mbox{and}\,\omega=\displaystyle\frac{16a}{9}-\displaystyle\frac{c^2}{4}.
\end{equation}
The computations are done for $-80\leq x,y\leq 100$ with different space and time steps as shown in the figures and tables corresponding. We fix also soliton parameters $a=0,01$, $c=0,1$ and the phase parameters $\varphi=\phi=0$. Therefore, the soliton pair $(u,v)$ is the exact solution of the inhomogeneous problem
\begin{equation}\label{InhpmogeneousContinuousProblem1}
\left\{\begin{array}{lll}
iu_{t}+\Delta u+|u|^{p-1}u+|v|^2u=G_1(x,y,t),\\
iv_{t}+\Delta v+|v|^{p-1}v+|u|^2v=G_2(x,y,t)
\end{array}
\right.
\end{equation}
where $G_1$ and $G_2$ are explicited respectively by
$$
G_1(x,y,t)=K^2\Bigl(sech_D(y,t)\Bigr)^{\frac{4}{3}}u(x,t)
$$
and
$$
G_2(x,y,t)=\biggl[K^2\Bigl(sech_B(x,t)\Bigr)^{\frac{8}{3}}+i\displaystyle\frac{8c\sqrt{a}}{3}T_D(y,t)\biggr]v(y,t).
$$
\section{Conclusion}
In this paper numerical study of a mixed cubic superlinear coupled Schrodinger system is considered. A non standard numerical scheme is developed to approximate the solutions in two-dimensional case by using a dynamical generalized Lyapunov-Sylvester algebraic operators. The discrete algebraic system is proved to be uniquely solvable, stable and fastly convergent. Numerical examples illustrating both phenomena of propagation and interaction of solitons are provided to show the efficiency of the numerical method.

\end{document}